\documentclass[a4,11pt]{article}

\usepackage{amsmath, amsfonts, amssymb, graphicx}
\usepackage{amsmath,mathtools,amssymb,amsfonts, amsthm}

\usepackage{times}
\usepackage{bm}
\usepackage{natbib}
\usepackage{color}

\usepackage[plain,noend]{algorithm2e}
\usepackage{url}

\newtheorem{theorem}{Theorem}[section]
\newtheorem{corollary}{Corollary}[section]

\makeatletter
\renewcommand{\algocf@captiontext}[2]{#1\algocf@typo. \AlCapFnt{}#2} 
\def\@algocf@capt@plain{top}
\renewcommand{\algocf@makecaption}[2]{%
  \addtolength{\hsize}{\algomargin}%
  \sbox\@tempboxa{\algocf@captiontext{#1}{#2}}%
  \ifdim\wd\@tempboxa >\hsize
    \hskip .5\algomargin%
    \parbox[t]{\hsize}{\algocf@captiontext{#1}{#2}}
  \else%
    \global\@minipagefalse%
    \hbox to\hsize{\box\@tempboxa}
  \fi%
  \addtolength{\hsize}{-\algomargin}%
}
\makeatother

\newcommand{\vs}{\boldsymbol{s}}

\title{ On F-modeling based Empirical Bayes Estimation of Variances }

\author{Yeil Kwon\\ Department of Mathematics, University of Central Arkansas\\ 201 Donaghey Avenue, Conway AR 72035, USA\\ Zhigen Zhao\\ Department of Statistics, Operations, and Data Science, Temple University\\ 1810 Liacouras Walk, Philadelphia, PA 19122, USA }

\begin{document}




\maketitle

\begin{abstract}
  We consider the problem of empirical Bayes estimation of multiple variances when provided with sample variances. Assuming an arbitrary prior on the variances, we derive different versions of the Bayes estimators using different loss functions. For one particular loss function, the resulting Bayes estimator relies on the marginal cumulative distribution function of the sample variances only. When replacing it with the empirical distribution function, we obtain an empirical Bayes version called  F-modeling based empirical Bayes estimator of variances. We provide theoretical properties of this estimator and further demonstrate its advantages through extensive simulations and real data analysis. 
\end{abstract}

{\bf Keywords: uniform convergence, empirical distribution function, selective inference.}

\section{Introduction}\label{sec:intro}

The empirical Bayes approach was introduced as a compound decision procedure in \cite{Robbins:1951} and has been widely studied thereafter (\citealp{Robbins:1956, Dvoretzky;Kiefer;Wolfowitz:1956, Efron:Morris:1972, Efron:Morris:1973, Efron:Morris:1975, Laird:Louis:1987, Jiang:Zhang:2009, Gu:Koenker:2017}). 
This approach plays an important role in the kinds of data analysis conducted during gene expression experiments, which often involve a large number of parallel inference problems.

The core idea of the empirical Bayes approach is to estimate the prior distribution either directly or indirectly using the available data, wherein the final inference is based on the posterior distribution when using this estimated prior. \cite{Efron:2014} classified empirical Bayes approaches as pursing one of two strategies: (i) f-modeling, which is modeling on the data scale; and (ii) g-modeling, which is modeling on the parameter scale.
Under f-modeling, the resulting empirical Bayes rule usually depends on the prior indirectly via the marginal probability density function; under g-modeling, the prior distribution is estimated and then plugged into  the posterior calculation. It is further commented in that paper that the g-modeling approach has been widely used in theoretical investigations (\citealp{Laird:Louis:1987, Morris:1983a, Jiang:Zhang:2009}), whereas the f-modeling approaches are more prevalent in applications (\citealp{Robbins:1956, Brown:Greenshtein:2009, Efron:2011}). 

The simultaneous estimation of variances and the covariance matrix has a long history, dating back to \cite{James:Stein:1961}. \cite{haff1980empirical} provided a parametric empirical Bayes estimator of the covariance matrix by assuming an inv-Wishart prior distribution on the covariance matrix. \cite{efron1976multivariate} proposed an estimator to dominate the sample covariance. \cite{wild1980loss} considered simultaneous estimation of the variances under different loss functions. \cite{robbins1982estimating} discussed a parametric empirical Bayes methods for scale mixture of Gaussians. \cite{champion2003empirical} considered the shrinkage estimator of  variances based on the Kullback-Leibler distance.

Heteroskedasticity is prevalent in many applications, such as microarray experiments, rendering the simultaneous estimation of variances even more important.
There have been many attempts to estimate these parameters with different approaches (\citealp{Tusher:Tibshirani:Chu:2001, Lonnstedt:Speed:2002, Storey:Tibshirani:2003, Lin:Nadler:Attie:Yandell:2003, Tong:Wang:2007, Gu:Koenker:2017}). Among these, there are a few widely used parametric empirical Bayes estimators which are widely used. When assuming an inverse gamma prior, \cite{Smyth:2004} developed a parametric empirical Bayes estimator of the variances. \cite{Cui:Hwang:Qiu:Blades:Churchill:2005} approximated both the chi-square distribution and the inverse gamma prior by  log-normal random variables and derived the exponential Lindley-James-Stein estimator. \cite{Lu:Stephens:2016} assumed that the prior of the variances follows a mixture of inverse gamma distributions to derive a flexible empirical Bayes estimator. These parametric empirical Bayes methods have the advantage of providing the full posterior distribution of the variances for further inference such as constructing credible intervals and performing hypothesis testing.
\cite{Gu:Koenker:2017} took the g-modeling approach by estimating the probability density function of the prior distribution using non-parametric maximum likelihood estimator (\citealp{koenker2014convex, kiefer1956consistency}).

In this work, we assume an arbitrary prior distribution $g(\sigma^2)$ for the variances to produce a nonparametric empirical Bayes estimator. When assuming some commonly used loss functions, we derive empirical Bayes estimators for the variances by modeling on the data scale. For a particular loss function, the resulting Bayes estimator depends only on  the marginal cumulative distribution function of the sample variances, $F(s^2)$. To the best of the authors' knowledge, this is the first estimator for the variances which relies on the marginal cumulative distribution function rather than the marginal probability density function.  To differentiate our method from the terminology used in \cite{Efron:2014}, we call this estimator an F-modeling based estimator. The advantage of the F-modeling based estimator is that one can simply replace the marginal cumulative distribution function with the empirical distribution function to obtain the proposed empirical Bayes version, which we call F-modeling based empirical Bayes estimator for the variances. The computation of the proposed method is instantaneous without any tuning parameters. 

It is known that the empirical distribution function converges to the true distribution function uniformly (\citealp{Dvoretzky;Kiefer;Wolfowitz:1956}). 
As shown in Section \ref{sec:analytic}, the proposed empirical Bayes estimator converges to the Bayes version uniformly over a set $\mathcal{D}_\delta=(0, D_\delta)$ where $D_\delta$ is a large value and tends to infinity when $\delta$ goes to zero. We impose this condition for technical reasons so as to prevent the denominator of the Bayes estimator to be arbitrarily small. It causes little practical concern because most often one would be interested in parameters corresponding to the small and moderate sample variances which fall in $\mathcal{D}_\delta$. 
We have also derived the estimator of the variances for the post selection inference and finite Bayes inference (\citealp{Efron:2019}).

\section{Empirical Bayes Estimator  for Variances}\label{sec:npeb}
Let $\sigma^2_{[1:N]}=(\sigma_1^2,\sigma_2^2,\cdots, \sigma_N^2)$ be the parameters of interest and $s^2_{[1:N]}=(s_1^2,s_2^2,\cdots, s_N^2)$ be the corresponding sample variances. In this paper, we consider the following model, 
\begin{equation}
\begin{cases}\label{eqn:model}
s_i^2 |\sigma_i^2 &\overset{ind}{\sim} p\left(s_i^2|\sigma_i^2\right)\sim \sigma_i^2\frac{\chi_k^2}{k}, \\
\sigma_i^2 &\overset{iid}{\sim} g\left(\sigma^2_i\right).
\end{cases}
\end{equation}
Here, $\chi^2_k$ denotes the random variable which follows a chi-square distribution with $k$ degrees of freedom. We assume an arbitrary prior $g(\sigma^2_i)$ on the variances. When integrating the variance out, the marginal probability density function of the sample variances is $f(s_i^2)=\int_0^\infty p(s_i^2|\sigma^2_i)g(\sigma^2_i)d\sigma^2_i$. Let
\begin{eqnarray}
F(s^2_i)=\int_0^{s^2_i} f(s^2_i)ds^2_i 
\end{eqnarray}
be the corresponding marginal cumulative distribution function of $s_i^2$'s. 

To derive the Bayes rule $\hat{\sigma}^2_{[1:N]}=(\hat{\sigma}_1^2,\hat{\sigma}_2^2,\cdots, \hat{\sigma}_N^2)$, a loss function must be specified. \cite{sinha1985inadmissibility} summarized many commonly used loss functions as follows:
\begin{eqnarray}\label{eq:loss}
\begin{array}{l} 
L_0\left( \sigma^2_{[1:N]},\hat\sigma^2_{[1:N]} \right) =\sum_{i=1}^N\left(\sigma_i^2- \hat \sigma_i^2\right)^2, \nonumber\\
L_1\left( \sigma^2_{[1:N]},\hat\sigma^2_{[1:N]} \right) =\sum_{i=1}^N\left(\frac{\sigma_i^2}{\hat \sigma_i^2}-1\right)^2,\\
L_1'\left( \sigma^2_{[1:N]},\hat\sigma^2_{[1:N]} \right) =\sum_{i=1}^N\left(\frac{\hat{\sigma}_i^2}{\sigma_i^2}-1\right)^2,\nonumber \\
L_2\left(  \sigma^2_{[1:N]},\hat\sigma^2_{[1:N]} \right) =\sum_{i=1}^N\left(\frac{\hat{\sigma}_i^2}{\sigma_i^2}-\ln\frac{\hat{\sigma}_i^2}{\sigma_i^2}-1\right).\nonumber \\
\end{array}
\end{eqnarray}

The squared error loss function, $L_0(\cdot)$, is not scale-invariant. The other three loss functions are scale-invariant. The loss function $L_1'(\cdot)$ is used in \cite{selliah1964estimation, ghosh1987}. The loss function $L_1'(\cdot)$ is equivalent to using $L_1(\cdot)$ when estimating the precision parameters (\citealp{ghosh1987}). 
The loss function $L_1'(\cdot)$ by nature favors under-estimation because ``underestimation has only a finite penalty, while overestimation has an infinite penalty"(\citealp{Casella:Berger:2001}). This could lead to an estimator which works extremely poor when focusing on the parameter with the smallest sample variance. On the contrary, both the loss function $L_1(\cdot)$ and Stein's Loss function $L_2(\cdot)$ have an infinite penalty for the underestimation. In addition, the loss function $L_2(\cdot)$ is tied to the Kullback-Leibler divergence and the entropy loss (\citealp{ghosh1987, wild1980loss, haff1977minimax,haff1980empirical}). A potential drawback of the loss function $L_1(\cdot)$ is that it imposes a finite penalty on the overestimation. 

In this article, we derive empirical Bayes estimators with respect to the scale-invariant loss functions $L_1'(\cdot)$, $L_1(\cdot)$, and $L_2(\cdot)$ by modeling on the data scale. We start with the loss function $L_1'(\cdot)$ where $\hat{\sigma}_{B,[1:N]}^{'2}=(\hat{\sigma}_{1,B}^{'2}, \hat{\sigma}_{2,B}^{'2},\ldots, \hat{\sigma}_{N,B}^{'2})$ is the corresponding Bayes rule. 
\begin{theorem}\label{theorem:1prime}
	Assume Model (\ref{eqn:model}) and the loss function $L_1'(\cdot)$, then
	\begin{equation}\label{eqn:est:varprime}
	\hat \sigma_{i,B}^{'2}
	=\frac{k(k-2)s_i^2f(s_i^2)-2ks_i^4f'(s_i^2) }{4s_i^4f''(s_i^2)-4(k-2)s_i^2{f'(s_i^2)}+{k(k-2)}f(s_i^2)}.
	\end{equation}
\end{theorem} 

Formula (\ref{eqn:est:varprime}) could be viewed as generalizing Tweedie's formula (\citealp{Efron:2011}) to the simultaneous estimation of variances. It is seen that the estimator $\hat{\sigma}_{i,B}^{'2}$ depends on the marginal probability density function $f(s_i^2)$,  its first and second derivatives.  We can get an empirical Bayes version by replacing $f(s_i^2)$ and its derivatives with the corresponding estimators using the kernel density estimator (\citealp{Brown:Greenshtein:2009}), or Lindsey's method (\citealp{Efron:2010b, Efron:2019}). We call this method the f-modeling based empirical Bayes estimator for variances:
\begin{equation}\label{eqn:fEBV}
\hat{\sigma}_{i,f-EBV}^{'2} 
=\frac{k(k-2)s_i^2\widehat{f(s_i^2)}-{2ks_i^4\widehat{f'(s_i^2)}}}{4s_i^4\widehat{f''(s_i^2)}-4(k-2)s_i^2\widehat{f'(s_i^2)}+{k(k-2)}\widehat{f(s_i^2)}}.
\end{equation} 
Next, consider the Stein's loss $L_2(\cdot)$ and let $\hat{\sigma}^2_{Stein,[1:N]}=(\hat{\sigma}^2_{1,Stein}, \hat{\sigma}^2_{2,Stein},\cdots,\hat{\sigma}^2_{N,Stein})$ be the corresponding Bayes rule. Then we have the following theorem. 
\begin{theorem}\label{theorem:stein}
	Assume Model (\ref{eqn:model}) and Stein's loss function $L_2(\cdot)$, then
	\begin{equation}\label{eqn:febv}
	\hat{\sigma}_{i,Stein}^{2} = \left( \frac{k-2}{ks_i^2} - \frac{2}{k}\cdot\frac{f'(s_i^2)}{f(s_i^2)}\right)^{-1}.
	\end{equation}
\end{theorem} 
When replacing $f(s^2)$ and $f'(s^2)$ with the corresponding estimators, we have the following f-modeling based empirical Bayes estimator of the variances when assuming Stein's loss:
\begin{equation}\label{eqn:fEBVS}
\hat{\sigma}_{i,f-EBVS}^{2} = \left( \frac{k-2}{ks_i^2} - \frac{2}{k}\cdot\frac{\widehat{f'(s_i^2)}}{\widehat{f(s_i^2)}}\right)^{-1}.
\end{equation}
When assuming Stein's loss, the empirical Bayes estimator does not require the estimation of the second derivative of the marginal probability density function. However, it still relies on the marginal density function and its first order derivative. 
The non-parametric estimation of the density function and its derivatives is a challenging problem, not to mention that the estimation accuracy on the tail becomes even worse. Additionally, the commonly used approaches such as the kernel density estimation relies on the choice of tuning parameters, which are difficult to choose in practice. 

Next, we consider the loss function $L_1(\cdot)$ and the corresponding Bayes decision rule $\hat{\sigma}_{B, [1:N]}^2=(\hat{\sigma}_{1,B}^2, \hat{\sigma}_{2,B}^2,\ldots, \hat{\sigma}_{N,B}^2)$. We have the following theorem.
\begin{theorem}\label{theorem:1}
	Assume Model (\ref{eqn:model}) and the loss function $L_1(\cdot)$. If 
	\begin{equation*}
	\int_{0}^{\infty}\,(s^2)^{-(\frac{k}{2}-2)}\,dF(s^2) <\infty\,\,\,\,\,
	\textrm{and}\,\,\,\,\,
	\int_{0}^{\infty}\,(s^2)^{-(\frac{k}{2}-1)} \,dF(s^2) <\infty,
	\end{equation*}
	then
	\begin{equation}\label{eqn:est:var}
	\hat \sigma_{i,B}^{2}
	=\frac{k}{2}\left\{ \frac{\int_{s_i^2}^{\infty}\,(s^2)^{-(\frac{k}{2}-2)}\,dF(s^2)}  {\int_{s_i^2}^{\infty}\,(s^2)^{-(\frac{k}{2}-1)} \,dF(s^2)} - s_i^2 \right\}.
	\end{equation}
\end{theorem} 

According to Model (\ref{eqn:model}), we know that
\[
\int_0^\infty (s^2)^{-(k/2-j)}dF(s^2)=\int_0^\infty\int_0^\infty C_k \frac{(s^2)^{j-1}}{\left(\sigma^2\right)^{k/2}} \exp\left(-\frac{ks^2}{2\sigma^2}\right)g(\sigma^2)d\sigma^2ds^2,\,\, j=1,2,
\]
where $C_k = \frac{ k^{k/2}}{\Gamma(k/2)2^{k/2}}$. When assuming an inverse gamma prior (\citealp{Smyth:2004}) and a mixture of inverse gamma prior (\citealp{Lu:Stephens:2016}), basic arithmetic calculations show that the conditions in the theorem hold.


Our F-modeling approach constructs a Bayes estimator of the variances which relies on $F(s^2)$, the cumulative distribution function of the sample variances. 
The advantage of using an F-modeling based estimator is that one can avoid the daunting task of estimating the marginal probability density function and its derivatives, which usually requires some kind of assumptions. Instead, to obtain an empirical Bayes version of the Bayes rule, we   simply replace $F(s^2)$ with the empirical distribution function $F_N(s^2) =\frac{1}{N}\sum_i I( s_i^2\le  s^2)$. After the substitution, we have the following proposed empirical Bayes estimator, which we refer to as the F-modeling based empirical Bayes estimator of the variances:
\begin{eqnarray}\label{est:5}
\hat \sigma^2_{i,F\textrm{-}EBV}
= \left\{\begin{array}{cl}
s_i^2, & \,\,\,\textrm{if $s_i^2 =\max\limits_{1\le j \le N} s_j^2$}, \\
\frac{k}{2}\left\{\frac{\sum_{s_j^2\ge s_i^2}(s_j^2)^{-(\frac{k}{2}-2)}}  {\sum_{s_j^2\ge s_i^2}(s_j^2)^{-(\frac{k}{2}-1)}} - s_i^2 \right\}, & \,\,\, \textrm{otherwise}.
\end{array}       \right.
\end{eqnarray}
The proposed estimator is calculated instantaneously and does not involve any tuning parameters. 


Return to Model (\ref{eqn:model}) with $g(\sigma^2)$ being arbitrary. Assume that one additional sample variance $s_0^2$ which is independent of $s_{[1:N]}^2$ has been observed. Let $\sigma_0^2$ be the corresponding variance which is assumed to be generated from $g(\sigma^2)$ and 
$ s_0^2 \sim \sigma_0^2\frac{\chi_k^2}{k}.$ The goal is to estimate $\sigma_0^2$ based on the posterior distribution $\sigma_0^2|s_0^2$. When $N$ goes to infinity, the prior distribution $g(\sigma^2)$ could be fully recovered and this reduces to the standard Bayes approach. For a finite $N$, this problem is called the finite Bayes inference (\citealp{Efron:2019}). Assume the loss function 
\begin{equation}\label{eqn:loss:FB}
L_1^{FB}(\hat{\sigma}_0^2, \sigma_0^2) = \left( \frac{\sigma_0^2}{\hat{\sigma}_0^2}-1\right)^2.
\end{equation}
Based on the proof of Theorem \ref{theorem:1}, we know that the Bayes rule is
\[
\hat{\sigma}^2_{0, B} =\frac{k}{2}\left\{ \frac{\int_{s_0^2}^{\infty}\,(s^2)^{-(\frac{k}{2}-2)}\,dF(s^2)}  {\int_{s_0^2}^{\infty}\,(s^2)^{-(\frac{k}{2}-1)} \,dF(s^2)} - s_0^2 \right\}.
\]
Consequently, we propose to estimate $\sigma_0^2$ by 
\begin{eqnarray}
\hat \sigma^2_{0,F-EBV}
= \left\{
\begin{array}{cl}
s_0^2, & \,\,\,\textrm{if $s_0^2\ge \max\limits_{1\le j\le N}s_j^2$}, \\
\frac{k}{2}\left\{\frac{\sum_{s_j^2\ge s_0^2}(s_j^2)^{-(\frac{k}{2}-2)}}  {\sum_{s_j^2\ge s_0^2}(s_j^2)^{-(\frac{k}{2}-1)}} - s_0^2 \right\}, & \,\,\, \textrm{otherwise}.
\end{array}
\right.
\end{eqnarray}

Similarly, we estimate $\sigma_0^2$ based on f-modeling methods by
\begin{equation}\label{eqn:FEBVS:finite}
\hat{\sigma}_{0,f-EBV}^{'2} 
=\frac{k(k-2)s_0^2\widehat{f(s_0^2)}-{2ks_0^4\widehat{f'(s_0^2)}}}{4s_0^4\widehat{f''(s_0^2)}-4(k-2)s_0^2\widehat{f'(s_0^2)}+{k(k-2)}\widehat{f(s_0^2)}}.
\end{equation}
and 
\begin{equation}\label{eqn:fEBVS:finite}
\hat{\sigma}_{0,f-EBVS}^{'2} = \left( \frac{k-2}{ks_0^2} - \frac{2}{k}\cdot\frac{\widehat{f'(s_0^2)}}{\widehat{f(s_0^2)}}\right)^{-1}.
\end{equation}

We can similarly construct estimators for variances relating to a set of indices, even if the indices have been chosen using the data.
Given the data $s^2_{[1:N]}=(s_1^2,s_2^2,\cdots, s_N^2)$, let $\mathcal{C}$ be the set of indices selected using a certain procedure. Our target is to estimate $\sigma_i^2, \forall i\in \mathcal{C}$ under the loss function
\begin{equation}\label{loss:selection}
L_1^\mathcal{PS}(\hat{\sigma}^2,\sigma^2) = \sum_{i\in\mathcal{C}} \left( \frac{\sigma_i^2}{\hat{\sigma}_i^2}-1 \right)^2.
\end{equation}
As an example, we might be interested in the variances corresponding to the $K$ smallest sample variances. In other words, order the sample variances $s_i^2$'s increasingly as $s_{(1)}^2 \le s_{(2)}^2\le \cdots \le s_{(N)}^2$. Let $\sigma^2_{(i)}$ be the parameter corresponding to $s_{(i)}^2$. Set $\mathcal{C}=\{i: s_i^2 \le s_{(K)}^2\}$. 

For any $i\in \mathcal{C}$, 
\[
\pi(\sigma_i^2|s^2_{[1:N]}, i\in\mathcal{C}) = \pi(\sigma_i^2|s^2_{[1:N]}).
\]
This implies that the posterior distribution of $\sigma_i^2$ when conditioning on both the data and the selection set is the same as the posterior distribution of $\sigma^2$ conditioning on the data. 
Consequently, the Bayes rule based on the selection remains the same and it is immune to the selection (\citealp{Dawid:1994}). We therefore propose to estimate $\sigma_i^2, i\in\mathcal{C}$ according to (\ref{est:5}) without adjustment. We would like to point out that this argument is true because the full data set is available for the post-selection inference. 
Otherwise, the Bayes rule might be affected by the selection. For instance, if only the data post the selection is available for further inference, then the Bayes rule needs to be corrected for such a selection rule. See \cite{Yekutieli:2012} for a full discussion on this issue.

\section{Theoretical Properties}\label{sec:analytic}

In this section, we study the theoretical properties of the proposed method. To ease our notation, we define two functions 
$l_1(s^2,u) = (s^2)^{-(k/2-2)}\mathbb{I}( s^2\geq u)$ and $l_2(s^2,u) = (s^2)^{-(k/2-1)}\mathbb{I}( s^2\geq u)$ where $\mathbb{I}(\cdot)$ is an indicator function. Then the Bayes decision rule and the proposed method can be respectively written as
\[\hat{\sigma}^2_{i,B} = \frac{k}{2}\left\{\frac{\int_0^\infty l_1(s^2,s_i^2)dF(s^2)}{\int_0^\infty l_2(s^2,s_i^2)dF(s^2)}-s_i^2\right\},
\,\,\,\, \textrm{and } \,\,\,\, \hat{\sigma}^2_{i,F-EBV} = \frac{k}{2}\left\{\frac{\int_0^\infty l_1(s^2,s_i^2)dF_N(s^2)}{\int_0^\infty l_2(s^2,s_i^2)dF_N(s^2)}-s_i^2\right\}.
\]
First, we study the numerator and denominator separately.

\begin{theorem}\label{theorem:2}
	Assume the same conditions in Theorem \ref{theorem:1} and $F(s^2)$ is continuous with the support of $(0,\infty)$, then
	\begin{equation*}
	\sup_{u} \left|\int_0^\infty l_1(s^2,u) dF_N(s^2)-\int_0^\infty l_1(s^2,u) dF(s^2) \right| \,\,\,\, \overset{a.s.}{\to} \,\,\,\,  0 ,
	\end{equation*}
	and
	\begin{equation*}
	\sup_{u} \left|\int_0^\infty l_2(s^2,u) dF_N(s^2)-\int_0^\infty l_2(s^2,u) dF(s^2) \right| \,\,\,\, \overset{a.s.}{\to} \,\,\,\,  0.
	\end{equation*}
\end{theorem}

This theorem implies that both the numerator and the denominator of the proposed empirical Bayes estimator converge to those of the Bayes rule uniformly.
However, it does not guarantee that the ratio converges uniformly. The reason is that the denominator $\int_0^\infty l_2(s^2, u)dF(s^2)$ converges to zero when $u$ goes to $\infty$. 
To prove that the proposed method converges to the Bayes estimator uniformly, we consider the set such that the denominator of the Bayes rule is greater than some positive number. Namely, for a number $\delta>0$, let $\mathcal{D}^\delta$ be a set defined as
\begin{equation}\label{eq:D1}
\mathcal{D}^\delta\equiv\left\{\,u\,\bigg|\,\int_{u}^\infty (s^2)^{-\left(\frac{k}{2}-1\right)} \,dF(s^2)>\delta\right\}.
\end{equation}
Since $\int_0^\infty (s^2)^{-(\frac{k}{2}-1)}dF(s^2) <\infty$, then $\mathcal{D}^\delta= (0, D_\delta)$ for some positive number $D_\delta$. We then have the following theorem:

\begin{theorem}\label{theorem:3}
	Assume the same conditions in Theorem \ref{theorem:2}, then
	\begin{equation*}	
	\sup_{s_i^2\in \mathcal{D}^\delta}\bigg| \hat\sigma^2_{i,F\textrm{-}EBV} -\hat\sigma^2_{i,B}\bigg| \,\,\,\,\overset{a.s.}{\to}\,\,\,\,0. 
	\end{equation*} 
\end{theorem}

The constant $D_\delta$ is a quantity depending on the marginal distribution function of the sample variances only and $D_\delta$ tends to infinity when $\delta$ tends to $0$. For any $0<\tau<1$, let $s^2_{[1:N]}$ be a random sample consisting of $N$ sample variances. Let $s_{\tau}^2$ be the $\tau$-th sample quantile. We can always choose $\delta$ sufficiently small, such that $\{s_i^2, s_i^2\le s_\tau^2\}\in \mathcal{D}^\delta$ with large probability. For a sample variance which doesn't fall in $\mathcal{D}^\delta$, one could estimate the corresponding parameter by this sample variances. Namely, we could modify the proposed estimator as
\begin{eqnarray}
\hat \sigma^{2}_{i,mF-EBV}
= \left\{
\begin{array}{cl}
s_i^2, & \,\,\,\textrm{if $s_i^2\ge s_{(\lceil N\tau\rceil)}^2$}, \\
\frac{k}{2}\left\{\frac{\sum_{s_j^2\ge s_i^2}(s_j^2)^{-(\frac{k}{2}-2)}}  {\sum_{s_j^2\ge s_i^2}(s_j^2)^{-(\frac{k}{2}-1)}} - s_i^2 \right\}, & \,\,\, \textrm{otherwise}.
\end{array}
\right.
\end{eqnarray}
In practice, especially when focusing on parameters with small sample variances, this modification does not make much difference. 

We can extend the result to the post-selection inference and finite Bayes inference.

\begin{corollary}\label{corrolary:1}
	Assume the same conditions in Theorem \ref{theorem:2}, then
	\begin{equation*}	
	\sup_{s_i^2\in \mathcal{D}^\delta, i\in\mathcal{C}}\bigg| \hat\sigma^2_{i,F\textrm{-}EBV} -\hat\sigma^2_{i,B}\bigg| \,\,\,\,\overset{a.s.}{\to}\,\,\,\,0.
	\end{equation*} 
\end{corollary}
As commented in Section \ref{sec:npeb}, the Bayes estimator is immune to the selection rule $\mathcal{C}$, and the empirical Bayes estimator could be a good approximation of the Bayes estimator. However, the discrepancy between these two widens when focusing on the selected case (\citealp{Pan:Huang:Hwang:2017}), and some correction is needed (\citealp{Hwang:Zhao:2013}). 
On the other hand, Corollary \ref{corrolary:1} indicates that the proposed F-modeling based empirical Bayes estimator converges to the corresponding Bayes version if $s_{i}^2\in \mathcal{D}^\delta, i\in\mathcal{C}$. 
In other words, we don't need to make further correction for the selection. 

Similarly, when considering the finite Bayes inference, the uniform convergence of the proposed estimator guarantees a good estimation as long as $s_0^2\in \mathcal{D}^\delta$.
\begin{corollary}\label{corrolary:2}
	Assume  the conditions in Theorem \ref{theorem:2}, then
	\begin{equation*}	
	\sup_{{s_0}^2\in \mathcal{D}^\delta}\bigg| \hat\sigma^2_{0,F\textrm{-}EBV} -\hat\sigma^2_{0,B}\bigg| \,\,\,\,\overset{a.s.}{\to}\,\,\,\,0.
	\end{equation*} 
\end{corollary}

\section{Numerical studies}\label{sec:application}

In this section, we compare the numerical performances of the proposed methods with existing methods, including the sample variance ($s^2$), exponential Lindley-James-Stein estimator (ELJS, \citealp{Cui:Hwang:Qiu:Blades:Churchill:2005}), Tong and Wang's method (TW, \citealp{Tong:Wang:2007}), Smyth method(\citealp{Smyth:2004}), variance adaptive shrinkage method (Vash, \citealp{Lu:Stephens:2016}), and  REBayes method (\citealp{Gu:Koenker:2017}). As suggested by a referee, we consider two more estimators based on the Smyth method and variance adaptive shrinkage method by considering the loss function $L_1(\cdot)$. Assume that the prior distribution $g(\sigma_i^2)$ in Model (\ref{eqn:model}) is inverse gamma $(a_0, b_0)$, then the posterior distribution of $\sigma_i^2$ is inverse gamma $(a_1,b_1)$ where $a_1=a_0+k/2$, $b_1=b_0+ks_i^2/2$. The hyper parameters $a_0$ and $b_0$ are estimated by using the method of moments (\cite{Smyth:2004}). The Smyth method, which minimizes $EL'_1(\cdot)$, is given as $\frac{b_1}{a_1}$. The modified Smyth method, which minimizes $EL_1(\cdot)$, is given as 
\[
\sigma_{i,mSmyth}^2 = \frac{E\left(\sigma_i^4|s_i^2\right)}{E\left(\sigma_i^2|s_i^2\right)}=\frac{b_1}{a_1-2}.
\]
Similarly, we include two versions of variance adaptive shrinkage estimators, the original version (Vash) and modified version (mVash) in our simulation studies.

Let $(\sigma_i^2, s_i^2), i=1,2,\cdots, N$ be the parameters and the sample variances be generated according to Model (\ref{eqn:model}) where the degrees of freedom $k$ is chosen as 5 and the prior $g(\sigma^2)$ is chosen from
\begin{enumerate}
	\item[] Setting I: $\sigma_i^2\sim$ inverse gamma distribution: $IG(a,1)$ where $a=10$ and $6$; 
	\item[] Setting II: $\sigma_i^2\sim$ Mixture of inverse gamma distributions: $0.2IG(a,1)+0.4IG(8,6) + 0.4IG(9, 19)$, where $a=10$ and $6$;
	\item[] Setting III: $\sigma_i^2=a$ with 0.4 probability and $1/a$ with 0.6 probability, where $a=3$ and $4$;
	\item[] Setting IV: $\sigma_i^2\sim$ Mixture of inverse Gaussian distributions: $0.4InvGauss(1/a, 1) + 0.6 InvGauss(a, a^4)$, where $a=2$ and $3$.
\end{enumerate}

For all simulations, we set $N=1,000$ and the number of replications as 500. For each replication, we generate the data $(\sigma_i^2,s_i^2)$ and order them according to the sample variances  increasingly. We consider three different selection rules: (i) the parameters corresponding to the 1\% smallest sample variances; (ii) the parameters corresponding to the 5\% smallest sample variances; and (iii) all the parameters.  
We calculate the estimated values based on the aforementioned methods. The risks associated with the loss function (\ref{loss:selection}) are calculated and reported in Table \ref{tab:Varcomp} and the table in the Appendix B. In our numerical studies, it is shown that two f-modelling estimators defined in (\ref{eqn:fEBV}) and (\ref{eqn:fEBVS}) perform poorly, and the results are not reported in the tables. 
The proposed F-modeling based empirical Bayes estimator performs the best among all the estimators considered. The modified Smyth method and modified variance adaptive shrinkage method perform similarly under these settings. Under Setting I when the prior of the variance is an inverse gamma distribution, the proposed method, the modified Smyth method and modified variance adaptive shrinkage method are essentially the same. 
However, for Settings II to IV when the prior distribution is not an inverse gamma distribution, 
the proposed method outperforms all other competing methods, including the modified Smyth method and the modified variance adaptive shrinkage method.

\begin{table}[h!]
	\centering
	\small{
		\begin{tabular}{|c|c|c|ccccccccc|}\hline	
			Setting &$a$          &	\%  & $s^2$&\textit{ELJS}&\textit{TW}&\textit{Smyth}&\textit{mSmyth}&\textit{Vash}&\textit{mVash}& \textit{REBayes}&\textit{Proposed}\\\hline
			&	&	1\%	&	2.60 & -0.48 &  -0.72 &  -0.90 & -1.06 & -0.87 & -1.06 & -0.65 & -1.06 \\
			I &10	&   5\% & 2.00 & -0.70 & -0.87 & -0.89 & -1.05 & -0.88 & -1.05 & -0.92 & -1.05\\	
			&	&	all	&  0.77 & -0.94 & -0.98 & -0.91 & -1.05 & -0.92 & -1.05 & -0.97 & -1.03\\	\hline
			&	&	1\%	&  2.34 & 1.05 & 0.45 & -0.14 & -0.21 & 0.87 & -0.10 & -0.05 & -0.22\\	
			II &10	&	5\%	&	1.79 & 0.62 & 0.17 & -0.10 & -0.20 & 0.74 & -0.11 & -0.06 & -0.22\\	
			&	&	all	&  0.75 & 0.01 & 0.00  & 0.14 & -0.43&  0.26&  -0.48 & -0.38 & -0.52\\	\hline
			&	&	1\%	&	 2.22 & 1.15 & 0.88 & -0.28 & -0.48 & -0.26 & -0.49 & -0.50 & -0.60\\	
			III &4	&	5\%	&	 1.72 & 0.74 & 0.53 & -0.06 & -0.36 & -0.05 & -0.37 & -0.22&  -0.39\\	
			&	&	all	&	 0.69 & 0.10 & 0.16 &  0.26 & -0.35 &  0.27 & -0.35 & -0.32 & -0.58\\	\hline
			&	&	1\%	&  2.28 & 1.26 & 0.97 & -0.08 & -0.28 & -0.06 & -0.28 & -0.13 & -0.28\\	
			IV & 4	&	5\%	&  1.73 & 0.77 & 0.53 & -0.13 & -0.28 & -0.11 & -0.29 & -0.22 & -0.32\\	
			&	&	all	& 0.72 & 0.14 & 0.20  & 0.29 & -0.34 &  0.30 & -0.34 & -0.30 & -0.56\\	
			\hline
	\end{tabular}}
	\caption{The $\log_{10}(risk)$ associated with the loss function (\ref{loss:selection}) of the different estimators for the variances under different simulation settings. For each setting, we consider three selection rule: (i) the parameters corresponding to the 1\% smallest sample variances; (ii) the parameters corresponding to the 5\% smallest sample variances; and (iii) all the parameters.
	}\label{tab:Varcomp}
\end{table}

Next, we consider the finite Bayes inference problem. Namely, for each generated data set $s^2_{[1:N]}$ and a new observation $s_0^2$, we calculate the estimated values based on different approaches and calculate the risk according to the loss function (\ref{eqn:loss:FB}). The risks are reported in Table \ref{tab:finiteBayes:1} and the table in Appendix B.  Overall, the proposed F-modeling based empirical Bayes estimator performs the best among all the estimators considered. The modified Smyth method and modified variance adaptive shrinkage method are essentially the same. 
Under Setting I when the prior of the variance is an inverse gamma distribution, the proposed method, the modified Smyth method and modified variance adaptive shrinkage method perform similarly with negligible differences. However, for Settings II to IV when the prior distribution is not an inverse gamma distribution, the proposed method 
outperforms all other competing methods.

\begin{table}[h!]
	\centering
	\small{
		\begin{tabular}{|c|c|ccccccccc|}
			\hline	
			Setting&$(a,b)$  & $s^2$&\textit{ELJS}&\textit{TW}&\textit{Smyth}&\textit{mSmyth}&\textit{Vash}&\textit{mVash}&\textit{REBayes}&\textit{Proposed}\\\hline
			I & 10	&  0.38 & 0.16 & -1.05 & -0.96 & -1.06 & -0.96 & -1.07 & -1.02 & -1.03\\ \hline
			II  & 10	& 0.36 & 0.14 & -0.11 & 0.01 & -0.48 & -0.02 & -0.5 & -0.51 & -0.55 \\ \hline
			III &4	&  0.92 & 0.72 & 0.23 & 0.23 & -0.36 & 0.25 & -0.36 & -0.31 & -0.47\\ \hline
			IV& 4	&   0.7 & 0.49 & 0.25 & 0.37 & -0.3 & 0.38 & -0.29 & -0.1 & -0.51	\\ \hline
	\end{tabular}}
	\caption{The $\log_{10}(risk)$ associated with the loss function (\ref{eqn:loss:FB}) of the different estimators for the {\it finite Bayes} inference problem. 
	}
	\label{tab:finiteBayes:1}
\end{table}

\section{Real data Analysis}\label{sec:realdata}
In this section, we apply different variance estimators to two microarray dataset: colon cancer \citep{Alon:Barkai:Notterman:1999} and Leukemia data \citep{Golub:Slonim:Tamayo:1999}. The colon cancer data contains gene expressions of genes ($N$=2,000) for 22 patients and 40 normal people. The leukemia data includes the expressions of genes ($N$ = 7,128) extracted from 72 patients with two types of leukemia: Acute Lymphoblastic Leukemia (47 patients) and Acute Myeloid Leukemia (25 patients). 
For the Leukemia data set, we first randomly split the subjects into two subgroups such that both subgroups contain similar numbers of subjects from the Acute Lymphoblastic Leukemia patients and Acute Myeloid Leukemia patients. For each sub-group, we then constructed $1-\gamma$ ($\gamma=0.05$) confidence intervals for $\theta_i$, the mean parameter of the $i$-th gene,  following the work of \cite{Hwang:Qiu:Zhao:2009} by considering 
\begin{equation*}
CI_i=\hat{\theta}_i \pm \sqrt{\hat{M}_i \hat\sigma_{i}^2}\cdot \sqrt{ z_{\gamma/2}^2 - \log \hat{M}_i},
\hat\theta_i =\hat{M}_iX_i+ (1-\hat{M}_i) \bar{X}, \hat{M}_i=\hat\tau^2/(\hat\sigma_{i}^2+\hat\tau^2),
\end{equation*}
and \[
\hat\tau^2=\max\left\{ \frac{1}{N}\sum_{i=1}^N(X_i-\hat\mu)^2-\frac{1}{N} \sum_{i=1}^N\sigma_i^2,\,\,\, \tau_0^2\right\}.
\]

We declare the $i$-th gene, where $i=1,2,\cdots,N$, to be significant if the corresponding interval does not enclose zero. We do the same for the other sub-group. We call the decision of the $i$-th gene discordant if the interval based on the first subgroup does (does not) enclose zero while the interval based on the second subgroup does not (does) enclose zero. If a decision is discordant, this implies that a significant conclusion based on one subgroup cannot be replicated by the other. We repeat these steps 500 times to calculate the average proportions of discordant decisions. We perform the same calculation for the colon cancer data by splitting the patients group and normal people group.

\begin{table}
	\centering
	\begin{tabular}{|c|ccccccccc|}
		\hline	
		Data & $s^2$& \small\textit{ELJS}& \small\textit{TW}& \small\textit{Smyth}& \small\textit{mSmyth}&\small\textit{Vash}&\small\textit{mVash}& \small\textit{REBayes}& \small\textit{F-EBV}\\
		\hline
		\textit{Colon}&0.27 &0.20 &0.20 & 0.20 &0.17 &0.20 &0.17 &0.19 &0.17\\
		\textit{Lukemia}&0.23 &0.15 &0.15 &0.15 & 0.13 & 0.16 &0.13  &0.14 &0.13\\\hline
	\end{tabular}
	\caption{The average percentages of discordant decisions of different intervals when applied to the Colon Cancer data and Leukemia Data based on 500 replications.
	}
	\label{tab:leukemia}
\end{table}

\begin{figure}
	\centering
	\includegraphics[width=70mm,height=45mm]{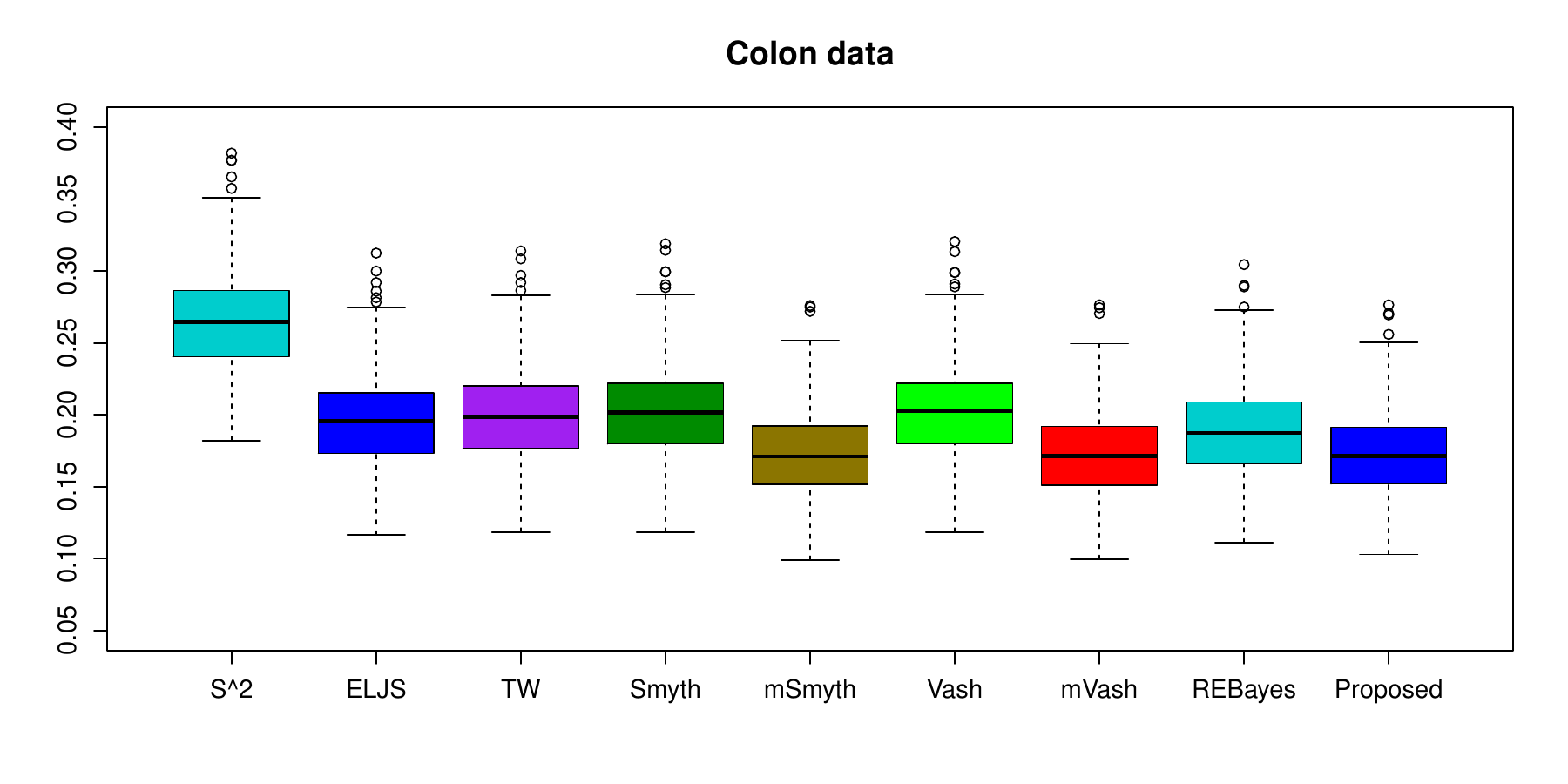}
	\includegraphics[width=70mm,height=45mm]{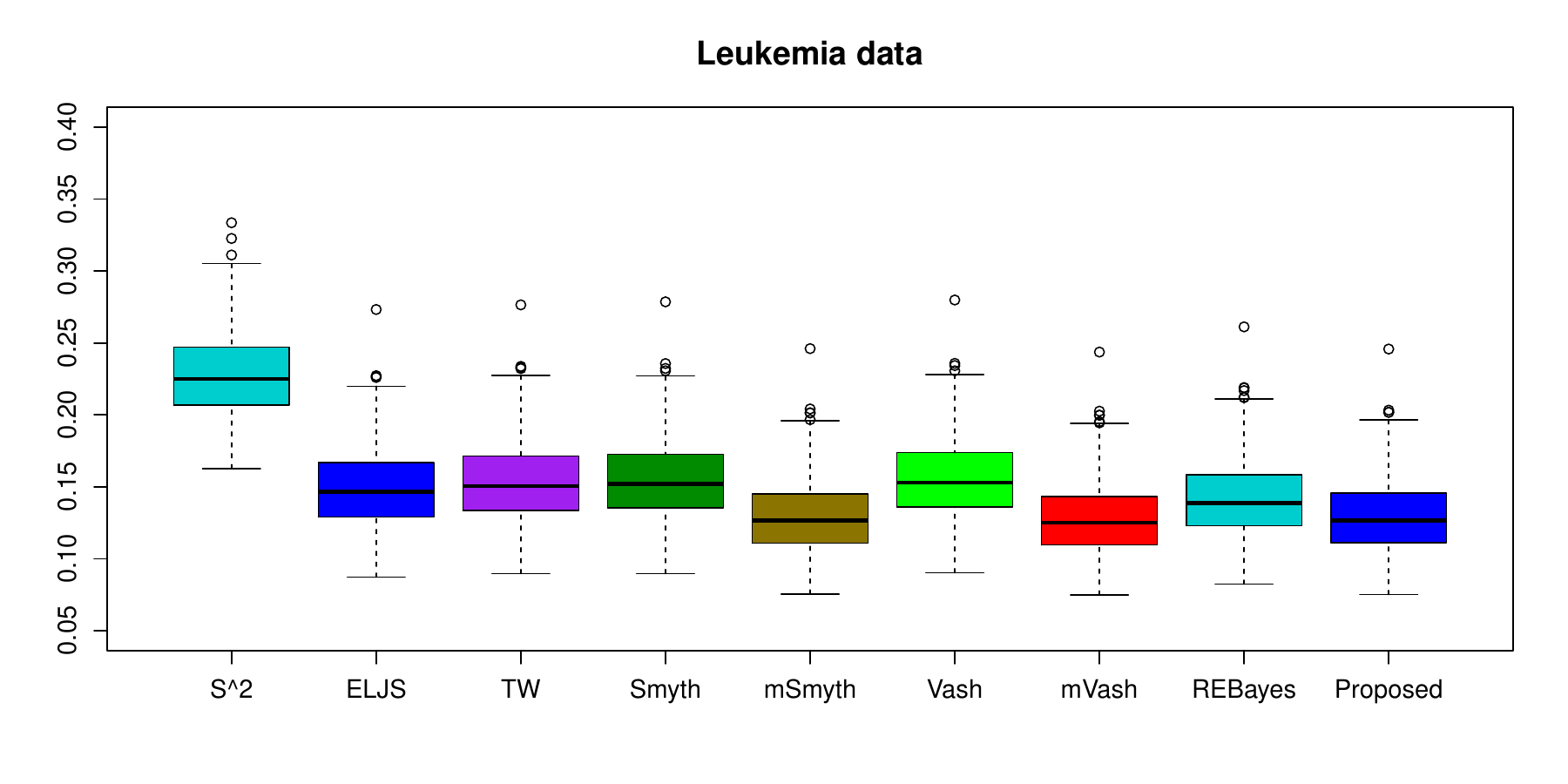} 
	\caption{The boxplots of the percentage of discordant decisions for colon cancer and leukemia data based on 500 replications. Left panel: Colon cancer data. Right panel: Leukemia data.}
	\label{fig:boxplot}
\end{figure}

In Figure \ref{fig:boxplot}, we plot the box-plots of the rate of discordant decisions. The average percentage of discordant decisions are reported in Table \ref{tab:leukemia}. It is seen that the proposed method, the modified Smyth and modified variance adaptive shrinkage estimator produce a similar number of discordance decisions. This number is substantially smaller than all the other competing methods. 

To further investigate why these three methods perform similarly, we test the hypothesis that the distribution of the sample variances is the convolution of a scaled chi-square distribution and an inverse gamma distribution. The Kolmogorov-Smirnov test statistics for the Colon data set and Leukemia data set are 0.014 and 0.017, respectively. The resulting p-values are 0.80 and  0.031, respectively. In other words, there is no evidence to reject the null hypothesis which states that the prior is an inverse gamma distribution for the colon data and there is only moderate evidence to reject the null hypothesis for the Leukemia data. It is expected to see similar performances for these three methods. 

The code for simulations and real data analysis are available on github (\url{https://github.com/zhaozhg81/FEBV}).

\section{Conclusion}\label{sec:conclusion}

The proposed method is developed under Model (\ref{eqn:model}) assuming a scaled chi-square distribution with equal degrees of freedom.  
The Bayes estimator in Theorem \ref{theorem:2} still applies when the degrees of freedom are different. However, the estimation of the cumulative distribution function requires that the sample variances are identically distributed. Therefore, the proposed method could not be directly applied to cases with unequal degrees of freedom. In practice, we take a slightly conservative approach by considering the smallest degrees of freedom as the common one. We would like to point out that many parametric empirical Bayesian approaches based on the g-modeling estimate the prior distribution explicitly and can handle unequal degrees of freedom.

In the real data analysis, we use the estimator of the variances as a plug-in estimator for inferring the mean parameters. One natural follow-up challenge to address is how to obtain a non-parametric empirical Bayes estimator of the means assuming arbitrary priors for both the means and the variances. Given the observed advantages of the F-modeling based approach, we would like to further extend this framework to broader settings in future research. We will further study the properties of the F-modeling based approach under the decision theoretical framework.

\section{Acknowledgement}
Zhigen Zhao's research is supported by the NSF grant IIS-1633283. Yeil Kwon's research is partially supported by the NSF grant IIS-1633283. The authors thank the AE and the reviewers for comments which helped substantially improve the quality of the paper. The authors also thank Mr. Matthew P MacNaughton for editing the manuscript.

\bibliographystyle{biometrika}
\bibliography{zhaozhg,References}

\begin{thebibliography}{43}
\expandafter\ifx\csname natexlab\endcsname\relax\def\natexlab#1{#1}\fi

\bibitem[{Alon et~al.(1999)Alon, Barkai, Notterman, Gish, Ybarra, Mack \&
  Levine}]{Alon:Barkai:Notterman:1999}
\textsc{Alon, U.}, \textsc{Barkai, N.}, \textsc{Notterman, D.~A.},
  \textsc{Gish, K.}, \textsc{Ybarra, S.}, \textsc{Mack, D.} \& \textsc{Levine,
  A.~J.} (1999).
\newblock Broad patterns of gene expression revealed by clustering analysis of
  tumor and normal colon tissues probed by oligonucleotide arrays.
\newblock In \textit{Proceedings of the National Academy of Science}, vol.~96.

\bibitem[{Brown \& Greenshtein(2009)}]{Brown:Greenshtein:2009}
\textsc{Brown, L.~D.} \& \textsc{Greenshtein, E.} (2009).
\newblock Nonparametric empirical {B}ayes and compound decision approaches to
  estimation of a high-dimensional vector of normal means.
\newblock \textit{The Annals of Statistics} \textbf{37}, 1685--1704.

\bibitem[{Casella \& Berger(2001)}]{Casella:Berger:2001}
\textsc{Casella, G.} \& \textsc{Berger, R.} (2001).
\newblock \textit{Statistical Inference}.
\newblock Duxbury Press, {S}econd ed.

\bibitem[{Champion(2003)}]{champion2003empirical}
\textsc{Champion, C.~J.} (2003).
\newblock Empirical {B}ayesian estimation of normal variances and covariances.
\newblock \textit{Journal of multivariate analysis} \textbf{87}, 60--79.

\bibitem[{Cui et~al.(2005)Cui, Hwang, Qiu, Blades \&
  Churchill}]{Cui:Hwang:Qiu:Blades:Churchill:2005}
\textsc{Cui, X.}, \textsc{Hwang, J.~T.}, \textsc{Qiu, J.}, \textsc{Blades,
  N.~J.} \& \textsc{Churchill, G.~A.} (2005).
\newblock Improved statistical tests for differential gene expression by
  shrinking variance components estimates.
\newblock \textit{Biostatistics} \textbf{6}, 59--75.

\bibitem[{Dawid(1994)}]{Dawid:1994}
\textsc{Dawid, A.~P.} (1994).
\newblock Selection paradoxes of {B}ayesian inference.
\newblock \textit{Institute of Mathematical Statistics Lecture Notes -
  Monograph Series} \textbf{24}, 211--220.

\bibitem[{Dvoretzky et~al.(1956)Dvoretzky, Kiefer \&
  Wolfowitz}]{Dvoretzky;Kiefer;Wolfowitz:1956}
\textsc{Dvoretzky, A.}, \textsc{Kiefer, J.} \& \textsc{Wolfowitz, J.} (1956).
\newblock Asymptotic minimax character of the sample distribution function and
  of the classical multinomial estimator.
\newblock \textit{The Annals of Mathematical Statistics} \textbf{27},
  642–669.

\bibitem[{Efron(2010)}]{Efron:2010b}
\textsc{Efron, B.} (2010).
\newblock \textit{Large-Scale Inference: Empirical Bayes Methods for
  Estimation, Testing, and Prediction}.
\newblock Cambridge Univ Pr.

\bibitem[{Efron(2011)}]{Efron:2011}
\textsc{Efron, B.} (2011).
\newblock Tweedie’s formula and selection bias.
\newblock \textit{Journal of the American Statistical Association}
  \textbf{106}, 1602--1614.

\bibitem[{Efron(2014)}]{Efron:2014}
\textsc{Efron, B.} (2014).
\newblock Two modeling strategies for empirical {B}ayes estimation.
\newblock \textit{Statistical Science} \textbf{29}, 285--301.

\bibitem[{Efron(2019)}]{Efron:2019}
\textsc{Efron, B.} (2019).
\newblock Bayes, oracle {B}ayes and empirical {B}ayes.
\newblock \textit{Statistical Science} \textbf{34}, 177--201.

\bibitem[{Efron et~al.(1976)Efron, Morris et~al.}]{efron1976multivariate}
\textsc{Efron, B.}, \textsc{Morris, C.} et~al. (1976).
\newblock Multivariate empirical {B}ayes and estimation of covariance matrices.
\newblock \textit{The Annals of Statistics} \textbf{4}, 22--32.

\bibitem[{Efron \& Morris(1972)}]{Efron:Morris:1972}
\textsc{Efron, B.} \& \textsc{Morris, C.~N.} (1972).
\newblock Limiting the risk of {B}ayes and empirical {B}ayes estimators. {II}.
  {T}he empirical {B}ayes case.
\newblock \textit{Journal of the American Statistical Association} \textbf{67},
  130--139.

\bibitem[{Efron \& Morris(1973)}]{Efron:Morris:1973}
\textsc{Efron, B.} \& \textsc{Morris, C.~N.} (1973).
\newblock Stein's estimation rule and its competitors---an empirical {B}ayes
  approach.
\newblock \textit{Journal of the American Statistical Association} \textbf{68},
  117--130.

\bibitem[{Efron \& Morris(1975)}]{Efron:Morris:1975}
\textsc{Efron, B.} \& \textsc{Morris, C.~N.} (1975).
\newblock Data analysis using {S}tein's estimator and its generalizations.
\newblock \textit{Journal of the American Statistical Association} \textbf{70},
  311--319.

\bibitem[{Ghosh \& Sinha(1987)}]{ghosh1987}
\textsc{Ghosh, M.} \& \textsc{Sinha, B.~K.} (1987).
\newblock Inadmissibility of the best equivariant estimators of the
  variance-covariance matrix, the precision matrix, and the generalized
  variance under entropy loss.
\newblock \textit{Statistics \& Risk Modeling} \textbf{5}, 201--228.

\bibitem[{Golub et~al.(1999)Golub, Slonim, Tamayo, Huard, Gaasenbeek, Mesirov,
  Coller, Loh, Downing, Caligiuri, Bloomfield \&
  Lander}]{Golub:Slonim:Tamayo:1999}
\textsc{Golub, T.~R.}, \textsc{Slonim, D.~K.}, \textsc{Tamayo, P.},
  \textsc{Huard, C.}, \textsc{Gaasenbeek, M.}, \textsc{Mesirov, J.~P.},
  \textsc{Coller, H.}, \textsc{Loh, M.~L.}, \textsc{Downing, J.~R.},
  \textsc{Caligiuri, M.~A.}, \textsc{Bloomfield, C.~D.} \& \textsc{Lander,
  E.~S.} (1999).
\newblock Molecular classification of cancer: class discovery and class
  prediction by gene expression monitoring.
\newblock \textit{Science} \textbf{286}, 531--537.

\bibitem[{Haff(1977)}]{haff1977minimax}
\textsc{Haff, L.} (1977).
\newblock Minimax estimators for a multinormal precision matrix.
\newblock \textit{Journal of Multivariate Analysis} \textbf{7}, 374--385.

\bibitem[{Haff(1980)}]{haff1980empirical}
\textsc{Haff, L.} (1980).
\newblock Empirical {B}ayes estimation of the multivariate normal covariance
  matrix.
\newblock \textit{The Annals of Statistics} \textbf{8}, 586--597.

\bibitem[{Hwang et~al.(2009)Hwang, Qiu \& Zhao}]{Hwang:Qiu:Zhao:2009}
\textsc{Hwang, J.~T.}, \textsc{Qiu, J.} \& \textsc{Zhao, Z.} (2009).
\newblock Empirical {B}ayes confidence intervals shrinking both means and
  variances.
\newblock \textit{Journal of the Royal Statistical Society. Series B}
  \textbf{71}, 265--285.

\bibitem[{Hwang \& Zhao(2013)}]{Hwang:Zhao:2013}
\textsc{Hwang, J.~T.} \& \textsc{Zhao, Z.} (2013).
\newblock Empirical {B}ayes confidence intervals for selected parameters in
  high dimension with application to microarray data analysis.
\newblock \textit{Journal of the American Statistical Association}
  \textbf{108}, 607--618.

\bibitem[{James \& Stein(1961)}]{James:Stein:1961}
\textsc{James, W.} \& \textsc{Stein, C.} (1961).
\newblock Estimation with quadratic loss.
\newblock In \textit{Proceedings of the Fourth Berkeley Symposium on
  Mathematical Statistics and Probability}, vol.~4. Berkeley, Calif.: Univ.
  California Press.

\bibitem[{Jiang \& Zhang(2009)}]{Jiang:Zhang:2009}
\textsc{Jiang, W.} \& \textsc{Zhang, C.~H.} (2009).
\newblock General maximum likelihood empirical {B}ayes estimation of normal
  means.
\newblock \textit{The Annals of Statistics} \textbf{37}, 1647--1684.

\bibitem[{Kiefer \& Wolfowitz(1956)}]{kiefer1956consistency}
\textsc{Kiefer, J.} \& \textsc{Wolfowitz, J.} (1956).
\newblock Consistency of the maximum likelihood estimator in the presence of
  infinitely many incidental parameters.
\newblock \textit{The Annals of Mathematical Statistics} \textbf{27}, 887--906.

\bibitem[{Koenker \& Gu(2017)}]{Gu:Koenker:2017}
\textsc{Koenker, R.} \& \textsc{Gu, J.} (2017).
\newblock {REB}ayes: An {R} package for empirical {B}ayes mixture methods.
\newblock \textit{Journal of Statistical Software} \textbf{82}, 1--26.

\bibitem[{Koenker \& Mizera(2014)}]{koenker2014convex}
\textsc{Koenker, R.} \& \textsc{Mizera, I.} (2014).
\newblock Convex optimization, shape constraints, compound decisions, and
  empirical {B}ayes rules.
\newblock \textit{Journal of the American Statistical Association}
  \textbf{109}, 674--685.

\bibitem[{Laird \& Louis(1987)}]{Laird:Louis:1987}
\textsc{Laird, N.~M.} \& \textsc{Louis, T.~A.} (1987).
\newblock Empirical {B}ayes confidence intervals based on bootstrap samples.
\newblock \textit{Journal of the American Statistical Association} \textbf{82},
  739--757.
\newblock With discussion and with a reply by the authors.

\bibitem[{Lin et~al.(2003)Lin, Nadler, Attie \&
  Yandell}]{Lin:Nadler:Attie:Yandell:2003}
\textsc{Lin, Y.}, \textsc{Nadler, S.~T.}, \textsc{Attie, A.~D.} \&
  \textsc{Yandell, B.~S.} (2003).
\newblock Adaptive gene picking with microarray data: detecting important low
  abundance signals.
\newblock \textit{The Analysis of Gene Expression Data: Methods and Software} ,
  291--312.

\bibitem[{L{\"o}nnstedt \& Speed(2002)}]{Lonnstedt:Speed:2002}
\textsc{L{\"o}nnstedt, I.} \& \textsc{Speed, T.} (2002).
\newblock Replicated microarray data.
\newblock \textit{Statistica Sinica} \textbf{12}, 31--46.
\newblock Special issue on bioinformatics.

\bibitem[{Lu \& Stephens(2016)}]{Lu:Stephens:2016}
\textsc{Lu, M.} \& \textsc{Stephens, M.} (2016).
\newblock Variance adaptive shrinkage (vash): flexible empirical bayes
  estimation of variances.
\newblock \textit{Bioinformatics} \textbf{32}, 3428--3434.

\bibitem[{Morris(1983)}]{Morris:1983a}
\textsc{Morris, C.~N.} (1983).
\newblock Parametric empirical {B}ayes confidence intervals.
\newblock In \textit{Proceedings of a Conference Conducted by the Mathematics
  Research Center, the University of Wisconsin–Madison}, Scientific
  inference, data analysis, and robustness. Orlando, FL: Academic Press.

\bibitem[{Pan et~al.(2017)Pan, Huang \& Hwang}]{Pan:Huang:Hwang:2017}
\textsc{Pan, J.}, \textsc{Huang, Y.} \& \textsc{Hwang, J.~G.} (2017).
\newblock Estimation of selected parameters.
\newblock \textit{Computational Statistics \& Data Analysis} \textbf{109},
  45--63.

\bibitem[{Robbins(1951)}]{Robbins:1951}
\textsc{Robbins, H.} (1951).
\newblock Asymptotically subminimax solutions of compound statistical decision
  problems.
\newblock In \textit{Proceedings of the {S}econd {B}erkeley {S}ymposium on
  {M}athematical {S}tatistics and {P}robability}, vol.~2. Berkeley and Los
  Angeles: University of California Press.

\bibitem[{Robbins(1956)}]{Robbins:1956}
\textsc{Robbins, H.} (1956).
\newblock An empirical {B}ayes approach to statistics.
\newblock In \textit{Proceedings of the {T}hird {B}erkeley {S}ymposium on
  {M}athematical {S}tatistics and {P}robability}, vol.~3. Berkeley and Los
  Angeles: University of California Press.

\bibitem[{Robbins(1982)}]{robbins1982estimating}
\textsc{Robbins, H.} (1982).
\newblock Estimating many variances.
\newblock In \textit{Statistical Decision Theory and Related Topics III}.
  Elsevier, pp. 251--261.

\bibitem[{Selliah(1964)}]{selliah1964estimation}
\textsc{Selliah, J.~B.} (1964).
\newblock \textit{Estimation and testing problems in a Wishart distribution}.
\newblock Ph.D. thesis, Department of Statistics, Stanford University.

\bibitem[{Sinha \& Ghosh(1985)}]{sinha1985inadmissibility}
\textsc{Sinha, B.~K.} \& \textsc{Ghosh, M.} (1985).
\newblock Inadmissibility of the best equivariant estimators of the
  variance-covariance matrix and the generalized variance under entropy loss.
\newblock Tech. rep., University of Pittsburgh.

\bibitem[{Smyth(2004)}]{Smyth:2004}
\textsc{Smyth, G.~K.} (2004).
\newblock Linear models and empirical {B}ayes methods for assessing
  differential expression in microarray experiments.
\newblock \textit{Statistical Applications in Genetics and Molecular Biology}
  \textbf{3}, Article 3.

\bibitem[{Storey \& Tibshirani(2003)}]{Storey:Tibshirani:2003}
\textsc{Storey, J.} \& \textsc{Tibshirani, R.} (2003).
\newblock {SAM} thresholding and false discovery rates for detecting
  differential gene expression in {DNA} microarrays.
\newblock In \textit{The Analysis of Gene Expression Data: Methods and
  Software}. New York: Springer, pp. 272--290.

\bibitem[{Tong \& Wang(2007)}]{Tong:Wang:2007}
\textsc{Tong, T.} \& \textsc{Wang, Y.} (2007).
\newblock Optimal shrinkage estimation of variances with applications to
  microarray data analysis.
\newblock \textit{Journal of the American Statistical Association}
  \textbf{102}, 113--122.

\bibitem[{Tusher et~al.(2001)Tusher, Tibshirani \&
  Chu}]{Tusher:Tibshirani:Chu:2001}
\textsc{Tusher, V.~G.}, \textsc{Tibshirani, R.} \& \textsc{Chu, G.} (2001).
\newblock Significance analysis of microarrays applied to the ionizing
  radiation response.
\newblock \textit{Proceedings of the National Academy of Science.} \textbf{98},
  5116--5121.

\bibitem[{Wild(1980)}]{wild1980loss}
\textsc{Wild, C.} (1980).
\newblock Loss functions and admissibility of normal variance estimators.
\newblock \textit{Canadian Journal of Statistics} \textbf{8}, 95--101.

\bibitem[{Yekutieli(2012)}]{Yekutieli:2012}
\textsc{Yekutieli, D.} (2012).
\newblock Adjusted {B}ayesian inference for selected parameters.
\newblock \textit{Journal of the Royal Statistical Society: Series B}
  \textbf{74}, 515--541.

\end{thebibliography}

\appendix

\section{Technical proofs.}

\vspace{3ex}
\textbf{Proof of Theorem \ref{theorem:1prime}.} According to the loss function $L_1'$,
\[
E\left[L_1'(\sigma^2_{[1:N]},\hat{\sigma}^2_{[1:N]})|s^2_{[1:N]}\right] =\sum_i \hat{\sigma}_i^4E\left[(\sigma_i^{2})^{-2}|s^2_{[1:N]}\right]-2\hat{\sigma}_i^2E\left[(\sigma_i^{2})^{-1}|s^2_{[1:N]}\right] + 1.
\]
Consequently, 
\[
\hat{\sigma}_{i,B}^{'2}=\frac{ E\left[(\sigma_i^{2})^{-1}|s^2_{[1:N]}\right]}{ E\left[(\sigma_i^{2})^{-2}|s^2_{[1:N]}\right] }.
\]
For ease of notation, we drop the subscript "$i$" in the proof. Recall that $p(s^2|\sigma^2)$ is the density function of $s^2|\sigma^2$ and $g(\sigma^2)$ is the prior distribution of $\sigma^2$. Note that 
\[
p(s^2|\sigma^2)= \frac{\left({s^2}\right)^{\frac{k}{2}-1}e^{-\frac{ks^2}{2\sigma^2}}}{\Gamma\left(\frac{k}{2}\right)2^\frac{k}{2}}\cdot \left(\frac{k}{\sigma^2}\right)^\frac{k}{2},\quad s^2>0.
\]
Then
\[
f(s^2) = \int p(s^2|\sigma^2)g(\sigma^2)d\sigma^2 = \int C_k \omega(s^2,\sigma^2)d\sigma^2,
\]
where 
\[
C_k = \frac{ k^{k/2}}{\Gamma(k/2)2^{k/2}},\quad\textrm{and}\quad \omega(s^2, \sigma^2) = \frac{(s^2)^{k/2-1}}{(\sigma^2)^{k/2}}\exp\left( -\frac{ks^2}{2\sigma^2}\right)g(\sigma^2).
\]
Take the derivative of $f(s^2)$ with respect to $s^2$, we know that
\begin{eqnarray*}
	&&f'(s^2)=\int C_k \frac{k-2}{2s^2}\omega(s^2,\sigma^2)d\sigma^2 - \frac{k}{2}\int C_k\frac{1}{\sigma^2}\omega(s^2,\sigma^2)d\sigma^2\\
	&=& \frac{k-2}{2s^2}f(s^2) - \frac{k}{2} E\left(\frac{1}{\sigma^2}\Big|s^2\right) \cdot f(s^2).
\end{eqnarray*}
This leads to 
\begin{equation}\label{eqn:1}
\frac{k}{2}E\left(\frac{1}{\sigma^2}\Big|s^2\right) \cdot f(s^2) = \frac{k-2}{2s^2}f(s^2) - f'(s^2).
\end{equation}
Take the second order derivative of $f(s^2)$ with respect to $s^2$, we have
\begin{eqnarray*}
	&&f''(s^2) = -\frac{k-2}{2s^4}f(s^2) + \frac{k-2}{2s^2}f'(s^2) -\frac{k}{2}\int C_k\frac{1}{\sigma^2}\left(\frac{k-2}{2s^2}-\frac{k}{2\sigma^2}\right)\omega(s^2,\sigma^2)d\sigma^2\\
	&=& -\frac{k-2}{2s^4}f(s^2) + \frac{k-2}{2s^2}f'(s^2) -\frac{k(k-2)}{4s^2}E\left(\frac{1}{\sigma^2}\Big|s^2\right)\cdot f(s^2) + \frac{k^2}{4}E\left(\frac{1}{\sigma^4}\Big|s^2\right)\cdot f(s^2).
\end{eqnarray*}
Consequently, 
\begin{equation}\label{eqn:2}
\frac{k^2}{4}E\left(\frac{1}{\sigma^4}\Big|s^2\right)\cdot f(s^2) = f''(s^2) -\frac{k-2}{s^2}f'(s^2) + \frac{k(k-2)}{4s^4}f(s^2).
\end{equation}
Combining (\ref{eqn:1}) and (\ref{eqn:2}), we know that 
\begin{eqnarray*}
	&&\sigma_{B}^{'2} = \frac{E\left[(\sigma^2)^{-1}|s^2\right]}{E\left[(\sigma^2)^{-2}|s^2\right]}=\frac{k(k-2)s^2f(s^2)-2ks^4f'(s^2)}{4s^4f''(s^2) -4(k-2)s^2f'(s^2) + k(k-2)f(s^2)}.
\end{eqnarray*}
\qed

\vspace{2ex}
\textbf{Proof of Theorem \ref{theorem:stein}.} For ease of notation, we drop the subscript ``$i$'' in the proof. Recall that Stein loss function is defined as
\begin{align*}
L_2(\sigma^2, \hat \sigma^2)=\frac{\hat\sigma^2}{\sigma^2}-\ln\left(\frac{\hat\sigma^2}{\sigma^2}\right)-1.
\end{align*}
Consequently,
\begin{align*}
EL_2(\sigma^2,\hat\sigma^2\big|{s^2})=\hat\sigma^2E\left[(\sigma^2)^{-1}|s^2\right]-\ln{\hat\sigma^2}+E(\ln{\sigma^2}|{s^2})-1.
\end{align*}
Therefore, the estimator $\hat{\sigma}_{Stein}^2$ which minimizes the above expression is 
\begin{align*}
\hat\sigma_{Stein}^2=\frac{1}{E\left[(\sigma^2)^{-1} |s^2\right]}
\end{align*}

According to the proof of Theorem \ref{theorem:1prime},
\[
\frac{k}{2}E\left[(\sigma^2)^{-1}|s^2\right]\cdot f(s^2) = \frac{k-2}{2s^2}f(s^2)-f'(s^2).
\]
Therefore, 
\[
\hat{\sigma}_{Stein}^2 = \frac{1}{E\left[(\sigma^2)^{-1}|s^2\right]} =  \left(\frac{k-2}{ks^2} -\frac{2f'(s^2)}{kf(s^2)}\right)^{-1}.
\]

\vspace{2ex}
\textbf{Proof of Theorem \ref{theorem:1}.} For ease of notation, we drop the subscript ``$i$'' in the proof. Recall that $p(s^2|\sigma^2)$ is the density function of $s^2|\sigma^2$ and $g(\sigma^2)$ is the prior distribution of $\sigma^2$. Note that $p(s^2|\sigma^2)$ is given as
\begin{equation}\label{eqn:def:sgivensigma}
p(s^2|\sigma^2)= \frac{\left({s^2}\right)^{\frac{k}{2}-1}e^{-\frac{ks^2}{2\sigma^2}}}{\Gamma\left(\frac{k}{2}\right)2^\frac{k}{2}}\cdot \left(\frac{k}{\sigma^2}\right)^\frac{k}{2},\quad s^2>0.
\end{equation}
Define $f(s^2)$, $n(s^2)$ and $h(s^2)$ as
\begin{equation}\label{eqn:def:m}
f(s^2)=\int_0^\infty \,p(s^2|\sigma^2)g(\sigma^2)\,d\sigma^2,
\end{equation}
\begin{equation}\label{eqn:def:n}
n(s^2)=\int_0^\infty \sigma^2 \,p(s^2|\sigma^2)g(\sigma^2)\,d\sigma^2,
\end{equation}
and     
\begin{equation}\label{eqn:def:h}
h(s^2)=\int_0^\infty (\sigma^2)^2 \,p(s^2|\sigma^2)g(\sigma^2)\,d\sigma^2.
\end{equation}
Note that $f(s^2)$ is the marginal distribution of $s^2$. Then
\begin{align*}
\hat\sigma_B^2=\frac{E\left[\,({\sigma^2})^2|s^2\,\right]}{E\left[\,{\sigma^2}|s^2\,\right]}
=\frac{\int_0^{\infty}(\sigma^2)^2 p(\sigma^2|s^2) \,d\sigma^2}{\int_0^{\infty}\sigma^2 \,p(\sigma^2|s^2) \,d\sigma^2}
=\frac{h(s^2)}{n(s^2)}.
\end{align*}
By differentiating $n(s^2)(s^2)^{-(\frac{k}{2}-1)}$ with respect to $s^2$, we have
\begin{align}\label{eq:deriv_n}
\left[n(s^2)(s^2)^{-(\frac{k}{2}-1)}\right]'&=-\frac{k}{2}f(s^2)(s^2)^{-(\frac{k}{2}-1)}.
\end{align}
Namely, 
\begin{equation}\label{eq:nsc1}		
{n(s^2)}{(s^2)^{-(\frac{k}{2}-1)}}=-\frac{k}{2}\int_0^{s^2}  f(t) \,t^{-(\frac{k}{2}-1)} dt + C,
\end{equation}
for some constant $C$.

On the other hand, from (\ref{eqn:def:n}), the left hand side of (\ref{eq:nsc1}) can be expressed as 
\begin{align}\label{eq:nsc2}	
n(s^2)(s^2)^{-({k}/{2}-1)}
&=\int_0^\infty {(s^2)^{-({k}/{2}-1)}\,\sigma^2 p(s^2|\sigma^2)g(\sigma^2)\,d\sigma^2}\\ \nonumber
&=\int_0^\infty \frac{\left(k/2\right)^{{k}/{2}}}{\Gamma\left(k/2\right)} {\left(\frac{1}{\sigma^2}\right)^{{k}/{2}-1} e^{-\frac{ks^2}{2\sigma^2}}}  \,g(\sigma^2)\,d\sigma^2.
\end{align}
From  (\ref{eq:nsc1}) and (\ref{eq:nsc2}),	as $s^2$ approaches to zero,
\begin{align*}
C=\lim_{s^2 \rightarrow 0} n(s^2)(s^2)^{-({k}/{2}-1)}=\frac{(k/2)^{k/2}}{\Gamma(k/2)}E\left(\frac{1}{\sigma^2}\right)^{k/2-1}=
\frac{k}{2}E\left(\frac{1}{S^2}\right)^{{k}/{2}-1},
\end{align*}
since, for $j=1,2$,
\begin{align*}
E\left(\frac{1}{S^2}\right)^{{k}/{2}-j}=\frac{{\left(k/2\right)}^{{k}/{2}-j}}{\Gamma\left({k}/{2}\right)}E\left(\frac{1}{\sigma^2}\right)^{{k}/{2}-j}.
\end{align*}	
Therefore, 
\begin{align*}		
{n(s^2)}{(s^2)^{-(\frac{k}{2}-1)}}	
&= -\frac{k}{2} \int_0^{s^2}f(t) \,t^{-(\frac{k}{2}-1)} dt + \frac{k}{2}\, E\left(\frac{1}{S^2}\right)^{\frac{k}{2}-1}\\
&=-\frac{k}{2}\int_0^{s^2}  f(t) \,t^{-(\frac{k}{2}-1)} dt + \frac{k}{2}\, \int_0^{\infty} f(t)\,t^{-\left(\frac{k}{2}-1\right)} \, dt\\
&=\frac{k}{2} \int_{s^2}^\infty \,t^{-(\frac{k}{2}-1)} dF(t).
\end{align*}
We can calculate $h(s^2)$ in the similar way. Take the first and second order derivatives of $h(s^2)(s^2)^{-(\frac{k}{2}-1)}$ with respect to $s^2$, we then have
\begin{equation}\label{eq:1}
\left[h(s^2)(s^2)^{-(\frac{k}{2}-1)}\right]'=-\frac{k}{2}n(s^2)(s^2)^{-(\frac{k}{2}-1)},
\end{equation}

\begin{equation}\label{eq:2nd:derivative}
\left[h(s^2)(s^2)^{-(\frac{k}{2}-1)}\right]''=\frac{k^2}{4}f(s^2)(s^2)^{-(\frac{k}{2}-1)}.
\end{equation}
Consequently,
\begin{align}\label{eq:3}		
\left[{h(s^2)}{(s^2)^{-(\frac{k}{2}-1)}}\right]'	
&=\int_0^{s^2} \frac{k^2}{4} f(t)t^{-(\frac{k}{2}-1)} \,dt + C_1,
\end{align}
and
\begin{align}\label{eq:4}		
{h(s^2)}{(s^2)^{-(\frac{k}{2}-1)}}	
&=\int_0^{s^2}\int_0^{y} \frac{k^2}{4}f(t) \,t^{-(\frac{k}{2}-1)} \,dt\,dy + C_1 s^2 +C_2 \nonumber\\
&=\frac{k^2}{4}\int_0^{s^2} f(t)\,t^{-(\frac{k}{2}-1)}(s^2-t)\,dt + C_1 s^2 +C_2,
\end{align}
for some constants $C_1$ and $C_2$.

From (\ref{eq:1}) and (\ref{eq:3}), as $s^2$ approaches to zero, similar argument shows that
\begin{align*}		
C_1=\lim_{s^2 \rightarrow 0} \left[{h(s^2)}{(s^2)^{-(\frac{k}{2}-1)}}\right]'	
=-\frac{k^2}{4}E\left(\frac{1}{S^2}\right)^{\frac{k}{2}-1}.
\end{align*}
Similarly, combine equations (\ref{eq:2nd:derivative}) and (\ref{eq:4}) and let $s^2$ approach to zero, 
\begin{align*}
C_2=\lim_{s^2 \rightarrow 0} h(s^2)(s^2)^{-(\frac{k}{2}-1)} = \frac{k^2}{4}E\left(\frac{1}{S^2}\right)^{\frac{k}{2}-2}.
\end{align*}
Thus, 
\begin{align*}		
\quad\,\, {h(s^2)}{(s^2)^{-(\frac{k}{2}-1)}}
&=\frac{k^2}{4}\int_0^{s^2} f(t)\,t^{-(\frac{k}{2}-1)}(s^2-t)\,dt 
-\frac{k^2}{4}E\left(\frac{1}{S^2}\right)^{\frac{k}{2}-1}s^2 + \frac{k^2}{4}E\left(\frac{1}{S^2}\right)^{\frac{k}{2}-2}\nonumber\\	
&=\frac{k^2}{4}\left[-\left(s^2\int_{s^2}^{\infty} f(t)\,t^{-\left(\frac{k}{2}-1\right)}\,dt\right)+ \left(\int_{s^2}^{\infty}f(t)\,t^{-(\frac{k}{2}-2)}\,dt\right) \right]\nonumber\\
&=\frac{k^2}{4}\int_{s^2}^{\infty}\,t^{-(\frac{k}{2}-1)}(t-s^2)\,dF(t). \nonumber
\end{align*}
Therefore,
\begin{align*}
\hat \sigma^2_B=\frac{h(s^2)}{n(s^2)}
=\frac{k}{2}\left[ \frac{\int_{s^2}^{\infty}\,t^{-(\frac{k}{2}-2)}\,dF(t)}  {\int_{s^2}^{\infty}\,t^{-(\frac{k}{2}-1)} \,dF(t)} - s^2 \right].
\end{align*}
\qed

\vspace{2ex}
\textbf{Proof of Theorem \ref{theorem:2}.}  We restate one of the momumental theorems in the empirical process, on which our proof is based (Blum, 1955; DeHardt, 1971).   

Let $\mathcal{F}$ be a set of measurable function. The bracket $[a,b]$ is the set of all the functions $l\in\mathcal{F}$ with $a\leq l\leq b$. An $\epsilon$-bracket is a bracket with $\Vert b-a \Vert\leq\epsilon$. The bracketing number $N_{[]}(\epsilon,\mathcal{F},L_1(P))$ is the minimum number of $\epsilon$-brackets with which $\mathcal{F}$ can be covered. 

{\it\textbf{Theorem} (Blum-DeHardt)\label{Blum}
	Let $\mathcal{F}$ be a class of measurable functions such that $ N_{[]}(\epsilon,\mathcal{F},L_1(P)) < \infty,\,\,$ for every $\epsilon>0$. Then $\mathcal{F}$ is P-Glivenko-Cantelli.
}

We only prove the part for the numerator and the denominator can be similarly done. 
Let $\mathcal{F}=\left\{l_1 : l_1(s^2,u)= (s^2)^{-(k/2-2)}\mathbb{I}( s^2>u),\,\, u >0\right\}$ and  $P l_1(s^2,u) =\int_{0}^{\infty} l_1(s^2,u)\,dF(s^2)=\int_{u}^{\infty} {s^2}^{-(k/2-2)}\,dF(s^2).$ It suffices to show that $\mathcal{F}$ is a \textit{P-Glivenko-Cantelli} class of functions. Since $F$ is continuous and $\int_{0}^\infty (s^2)^{-({k}/{2}-2)} \, dF(s^2) < \infty,\,\,$ for any $\,\,\epsilon>0,\,\,$ a collection of real numbers $0=v_0<v_1<v_2<\cdots <v_m = \infty$ can be found such that 
\begin{align*}
Pl_1(s^2,v_{j-1})-Pl_1(s^2,v_j)
&=\int_{v_{j-1}}^{\infty}  (s^2)^{-(k/2-2)}\, dF(s^2) -\int_{v_j}^{\infty} (s^2)^{-(k/2-2)}\, dF(s^2)\\
&=\int_{v_{j-1}}^{v_j} (s^2)^{-(k/2-2)}\, dF(s^2)\\
&\leq \epsilon 
\end{align*}
for all $1 \leq j \leq m$, with
\begin{equation*}
Pl_1(s^2,v_m^{-})=\lim_{v_m\uparrow {\infty}} Pl_1(s^2,v_m)=\lim_{v_m\uparrow {\infty}}\int_{v_m}^{\infty} (s^2)^{-(k/2-2)}\, dF(s^2)=0.
\end{equation*}

Consider the collection of brackets $\{[a_j,b_j],1 \leq j \leq m\}$, with $a_j(s^2) ={s^2}^{-(k/2-2)}\,\mathbb{I}(s^2 > v_{j})$ and $b_j(s^2)={s^2}^{-(k/2-2)}\,\mathbb{I}(s^2>v_{j-1})$. Now each $l_1 \in \mathcal{F}$ is in at least one bracket and $\vert a_j-b_j \vert_{P} = Pl_1(s^2,v_{j-1})-Pl_1(s^2,v_{j}^{-}) \leq \epsilon $ for all $1 \leq j \leq m$. 
Thus, by Blum-DeHardt theorem , $\mathcal{F}$ is a \textit{P-Glivenco-Cantelli Class} of functions. \qed\\

\vspace{2ex}

\textbf{Proof of Theorem \ref{theorem:3}.}  Let 
\begin{equation*}
A_N(s_i^2)=\int_0^\infty l_1(s^2,s_i^2) \,dF_N(s^2),\qquad
A(s_i^2)=\int_0^\infty  l_1(s^2,s_i^2) \,dF(s^2),
\end{equation*}
and
\begin{equation*}
B_N(s_i^2)=\int_0^\infty l_2(s^2,s_i^2)  \,dF_N(s^2),\qquad
B(s_i^2)=\int_0^\infty l_2(s^2,s_i^2)  \,dF(s^2).
\end{equation*}
According to the proof of Theorem \ref{theorem:2}, $\sup_{s_i^2\in R}|A_N(s_i^2)-A(s_i^2)|\to 0\,\,$ and $\sup_{s_i^2\in R}|B_N(s_i^2)-B(s_i^2)|\to 0\,\,a.s.$.
Let $L=\inf_{s_i^2\in D^\delta}\{ B(s_i^2)\}$. Then for any $\epsilon>0$, when $N$ is sufficiently large
\[
\inf_{s_i^2\in D^\delta}B_N(s_i^2) \ge L-\epsilon\,\,\,\,\, a.s.,
\]
and $ \sup_{s_i^2\in R}A_N(s_i^2) \le C,\,\,\, a.s.$ for some constant $C$. Then

\begin{eqnarray*}
	&&\sup_{s_i^2\in D^\delta}\left|\hat{\sigma}_{i,F\textrm{-}EBV}^2 - \hat{\sigma}_{i,B}\right| \\
	&=& \sup_{s_i^2\in D^\delta}\left|\,\frac{A_N(s_i^2)}{B_N(s_i^2)}-\frac{A(s_i^2)}{B(s_i^2)}\,\right|\\
	&=&\sup_{s_i^2\in D^\delta} \left|\, \frac{A_N(s_i^2)(B(s_i^2)-B_N(s_i^2))}{B_N(s_i^2)B(s_i^2)}+ \frac{A_N(s_i^2)-A(s_i^2)}{B(s_i^2)}\,\right|  \\
	&\le & \frac{C}{L^2}\sup_{s_i^2\in D^\delta}\left|B(s_i^2)-B_N(s_i^2)\right| + \frac{1}{L} \sup_{s_i^2\in D^\delta}\left |A(s_i^2)-A_N(s_i^2)\right| \to 0, \,\,\,a.s..
\end{eqnarray*}
\qed

\vspace{3ex}

\newpage
\section{Additional simulation results}\label{sec:additional:sim}
In this section, we include additional simulation results which are not listed in the paper due to the page limit. The numerical results consist of four parts: (a) results of variance estimation post-selection; and (b) results of Finite Bayes inference problem. 
\\

\noindent {\bf (a) Results of variance estimation post-selection.}\\
To help the readers, we restate the simulation settings here. Let $\sigma_i^2$'s be the parameters, and the sample variances $s_i^2$'s are generated according to Model \ref{eqn:model} where the degrees of freedom $k$ is chosen as 5. We consider the following different choices of the prior $g(\sigma^2)$:
\begin{itemize}
	\item[] Setting I: $\sigma_i^2\sim$ inverse gamma distribution: $IG(a,1)$ where $a=10$ and $6$; 
	\item[] Setting II: $\sigma_i^2\sim$ Mixture of inverse gamma distributions: $0.2IG(a,1)+0.4IG(8,6) + 0.4IG(9, 19)$, where $a=10$ and $6$;
	\item[] Setting III: $\sigma_i^2=a$ with 0.4 probability and $1/a$ with 0.6 probability, where $a=3$ and $4$;
	\item[] Setting IV: $\sigma_i^2\sim$ Mixture of inverse Gaussian distributions: $0.4InvGauss(1/a, 1) + 0.6 InvGauss(a, a^4)$, where $a=2$ and $3$.
\end{itemize}

After generating the data, order the sample variances increasingly. We consider three different selection rules: (i) select the parameters corresponding to the 1\% smallest sample variances; (ii) select the parameters corresponding to the 5\% smallest sample variances; and (iii) all the parameters. We report $\log_{10}(risk)$ in Table \ref{tab:Varcomp:2}.

\begin{table}
	\centering
	\small{
		\begin{tabular}{|c|c|c|ccccccccc|}\hline	
			Setting &$a$          &	\%  & $s^2$&\textit{ELJS}&\textit{TW}&\textit{Smyth}&\textit{mSmyth}&\textit{Vash}&\textit{mVash}& \textit{REBayes}&\textit{Proposed}\\\hline
			&	&	1\%	&	2.53 & -0.13 & -0.57 & -0.67 & -0.87 & -0.56 & -0.87 & -0.68 & -0.87\\
			I &6	&   5\% &  1.95 & -0.40 & -0.70 & -0.66 & -0.87 & -0.59 & -0.87 & -0.82 & -0.87\\	
			&	&	all	&    0.74 & -0.70 & -0.71 & -0.66 & -0.87 & -0.66 & -0.87 & -0.78 & -0.86\\	\hline
			&	&	1\%	&   2.44  & 1.01  & 0.41 & -0.06 & -0.24 & 0.81 & -0.17 & -0.07 & -0.22\\	
			II &6	&	5\%	&	1.88  & 0.57 &  0.10 & -0.08 & -0.26 & 0.64 & -0.19 & -0.17 & -0.26 \\	
			&	&	all	&  0.77 & -0.09 & -0.11 &  0.02 & -0.49 & 0.12 & -0.51 & -0.44 & -0.54 \\	\hline
			&	&	1\%	&	 2.33  & 1.02 &  0.57 & -0.14 & -0.47 & 0.72 & -0.22 & -0.36 & -0.49\\	
			III &3	&	5\%	&	1.78 &  0.57 &  0.20 & -0.14 & -0.42&  0.57 & -0.22&  -0.34 & -0.44 \\	
			&	&	all	&	0.70 & -0.05 & -0.04 &  0.10 & -0.44 & 0.17 & -0.47 & -0.44 & -0.61 \\	\hline
			&	&	1\%	&  2.32 & 1.06 & 0.55 & -0.22 & -0.29 & 0.69 & -0.12 & -0.27&  -0.34\\	
			IV & 3	&	5\%	&  1.77 & 0.62 & 0.23 & -0.14 & -0.28 & 0.61 & -0.13 & -0.24 & -0.31\\	
			&	&	all	& 0.73 & 0.02 & 0.03  & 0.17 & -0.40 & 0.25 & -0.45 & -0.38 & -0.56\\	
			\hline
	\end{tabular}}
	\caption{The $\log_{10}(risk)$ associated with the loss function (\ref{loss:selection}) of the different estimators for the variances under different simulation settings. For each setting, we consider three selection rule: (i) the parameters corresponding to the 1\% smallest sample variances; (ii) the parameters corresponding to the 5\% smallest sample variances; and (iii) all the parameters.
	}\label{tab:Varcomp:2}
\end{table}

\vspace{20pt}
\noindent {\bf (b) Results of finite Bayes inference problem.}\\

Next, we consider the {\it finite Bayes inference} problem. Namely, for each generated data set $\vs^2$ and a new observation $s_0^2$, we calculate the estimated values based on different approaches and calculate the loss according to the loss function (\ref{eqn:loss:FB}). We calculate the risk based on 500 replications and reported the results in Table \ref{tab:finiteBayes:2}.
\\

\begin{table}
	\centering
	\small{
		\begin{tabular}{|c|c|ccccccccc|}
			\hline
			Setting&$(a,b)$  & $s^2$&\textit{ELJS}&\textit{TW}&\textit{Smyth}&\textit{mSmyth}&\textit{Vash}&\textit{mVash}&\textit{REBayes}&\textit{Proposed}\\\hline
			I & 6	&    0.3 & 0.07&  -0.86 & -0.81 &   -1 & -0.8 &   -1 & -0.91&  -0.98\\ \hline
			II  & 6	&   0.64 & 0.43 & -0.18 & -0.04 & -0.53&  -0.02 & -0.54 & -0.52&  -0.59\\ \hline
			III &3	&   0.92 & 0.72 & -0.02 & 0.06 & -0.46 & 0.18 & -0.48 & -0.54 & -0.61\\ \hline
			IV& 3	&   0.43 & 0.21 & -0.08 & 0.08 & -0.43 & 0.16 & -0.46 & -0.44 & -0.59 \\ \hline
	\end{tabular}}
	\caption{The $\log_{10}(risk)$ associated with the loss function (\ref{eqn:loss:FB}) of the different estimators for the {\it finite Bayes} inference problem. 
	}
	\label{tab:finiteBayes:2}
\end{table}

\end{document}